\documentclass[12pt]{amsart}

\usepackage{amsmath,amssymb,amsthm, latexsym}
\usepackage{amsfonts}
\usepackage[mathscr]{eucal}
\usepackage{color}

\newtheorem{theorem}{Theorem}[section]
\newtheorem{lemma}[theorem]{Lemma}

\theoremstyle{definition}
\newtheorem{defi}[theorem]{Definition}

\theoremstyle{remark}
\newtheorem{remark}[theorem]{Remark}

\numberwithin{equation}{section}

\newcommand{\RR}[1]{\mathbb{#1}}

\newcommand{\rd}{{\mathbb R^d}}

\newcommand{\rr}{{\mathbb R}}

\def\R{{\mathbb R}}

\def\E{{\mathbb E}}
\def\P{{\mathbb P}}

\def\A{{\mathcal A}}
\def\L{{\mathcal L}}

\def\sE{{\mathcal E}}
\def\sF{{\mathcal F}}

\def\wh{\widehat}

\def\l{{\langle}}
\def\r{\rangle}
\def\<{{\langle}}
\def\>{\rangle}

\def\eps{\varepsilon}

\begin{document}

\title{\bf Space-time fractional diffusion on bounded domains}

\author{Zhen-Qing Chen}
\address{Zhen-Qing Chen, Department of Mathematics, University of Washington, Seattle, WA 98195, USA}
\email{zchen@math.washington.edu}

\author{Mark M. Meerschaert}
\address{Mark M. Meerschaert, Department of Statistics and Probability,
Michigan State University, East Lansing, MI 48823.}
\email{mcubed@stt.msu.edu}
\urladdr{http://www.stt.msu.edu/$\sim$mcubed/}
\thanks{Research of Zhen-Qing Chen is partially supported
by NSF Grants DMS-0906743 and DMR-1035196.
Research of Mark M. Meerschaert was partially
supported by NSF grants DMS-1025486, DMS-0803360, EAR-0823965 and NIH grant R01-EB012079-01.}

\author{Erkan Nane}
\address{Erkan Nane, Department of Mathematics and Statistics, 221 Parker Hall, Auburn University, Auburn, AL 36849.}
\email{ezn0001@auburn.edu}

\begin{abstract}
Fractional diffusion equations replace the integer-order derivatives in space and time by their fractional-order analogues.  They are used in physics to model anomalous diffusion.  This paper develops strong solutions of space-time fractional diffusion equations on bounded domains, as well as probabilistic representations of these solutions, which are useful for particle tracking codes.
\end{abstract}

\keywords{Fractional derivative; anomalous diffusion; probabilistic representation, strong solution; Cauchy problem; bounded domain}

\maketitle

\section{Introduction}

The traditional diffusion equation $\partial_t u=\Delta u$ describes
a cloud of spreading particles at the macroscopic level. The point source solution is a Gaussian probability density that predicts the relative particle concentration.  Brownian motion provides a microscopic picture,
describing the paths of individual particles.  A
Brownian motion, killed or stopped upon leaving a domain, can be used to solve
Dirichlet boundary value problems for the heat equation, as well as some elliptic equations \cite{bass,davies}.
The space-time fractional diffusion equation $\partial^\beta_t u=\Delta^{\alpha/2} u$
with $0<\beta<1$ and $0<\alpha<2$ models anomalous diffusion \cite{Zsolution}.   The fractional derivative in time can be used
to describe particle sticking and trapping phenomena.
The fractional space derivative models long particle jumps.  The combined effect produces a concentration profile with a sharper peak, and heavier tails.  This paper studies strong solutions, and probabilistic representations of solutions, for the space-time diffusion equation on bounded domains.  Our main result is Theorem \ref{thm:5.6}. Strong solutions are obtained by separation of variables, combining the Mittag-Leffler solution to the time-fractional problem with an eigenfunction expansion of the fractional Laplacian on bounded domains.   The probabilistic representation of solutions involves an inverse stable subordinator time change, resulting in a non-Markovian process.  Fractional diffusion equations are becoming popular in many areas of application \cite{GorenfloSurvey,MetzlerKlafter}.  In these applications, it is often important to consider boundary value problems.  Hence it is useful to develop solutions for space-time fractional diffusion equations on bounded domains with Dirichlet boundary conditions.

\section{Random walks and stable processes}

A random walk $S_t=Y_1+\cdots+Y_{[t]}$, a sum of independent and identically distributed $\R^d$-valued random vectors, is commonly used to model diffusion in statistical physics.  Here $[t]$ denotes the largest integer not exceeding $t$, and $S_n$ represents the location of a random particle at time $n$.  Suppose the distribution of $Y$ is spherically symmetric.  If $\sigma^2:=\E \left[  | Y_1 |^2 \right]$ is finite and
$\E [Y_1]=0$, Donsker's invariance principle  implies that as $\lambda\to \infty$, the random process $\{\lambda^{-1/2}S_{\lambda t}, \, t\geq 0\}$ converges weakly in the Skorohod space to a Brownian motion $\{B_t, \, t\geq 0\}$ with $\E [ B^2_1]=\sigma^2$.  If the step random variable $Y_1$ is spherically symmetric, and $\P( | Y_1 |>x)\sim C x^{-\alpha}$ as $x\to\infty$ for some $0<\alpha<2$ and $C>0$, then $\E \left[ |Y_1 |^2 \right]$ is infinite,
and the extended central limit theorem tells us that
$\{\lambda^{-1/\alpha}S_{\lambda t}, t\geq 0\}$ converges weakly to a rotationally symmetric $\alpha$-stable L\'evy motion $\{A_t, t\geq 0\}$
 with
$$
\E[e^{i \xi \cdot  A_t }]=e^{-C_0 |\xi |^\alpha t }
 \qquad \hbox{for every } \xi \in \R^d \hbox{ and } t\geq 0,
$$
 where the constant $C_0$ depends only on $C$ and the dimension $d$, see \cite{RVbook}.  A simple rescaling in space yields a standard stable process with $C_0=1$.  Since $\{\lambda^{1/\alpha} A_t, t\geq 0\}$ has the same
 distribution as $\{A_{\lambda t}, t\geq 0\}$,  stable L\'evy motion represents a model for anomalous super-diffusion, where particles spread faster than a Brownian motion \cite{fadePRE}.

If we impose a random waiting time $T_n$ before the $n$th random walk jump, then the position of the particle at time $T_n=J_1+\cdots+J_n$ is given by $S_n$.  The number of jumps by time $t>0$ is $N_t=\max\{n:T_n\leq t\}$, so the position of the particle at time $t>0$ is $S_{N_t}$, a subordinated process.  If $ \P(J_n>t)\sim C t^{-\beta}$ as $t\to\infty$ for some $0<\beta<1$, then the scaling limit of $c^{-1/\beta}T_{[ct]}\Rightarrow Z_t$ as $c\to \infty$ is a strictly increasing stable L\'evy motion with index $\beta$, sometimes called a stable subordinator.  The jump times $T_n$ and the number of jumps $N_t$ are inverses: $\{N_{t}\geq n\}=\{T_n\leq t\}$.
 \cite[Theorem 3.2]{limitCTRW} shows that
 $\{c^{-\beta}N_{ct}, t\geq 0\}$ converges weakly to the process
 $\{E_t, t\geq 0\}$, where $E_t=\inf\{x:Z_x> t\}$.
 In other words, the scaling limits are also inverses:
  $\{E_t\leq x\}=\{Z_x\geq t\}$.
 Now $N_{ct}\approx c^{\beta}E_t$, and \cite[Theorem 4.2]{limitCTRW}
 shows that the scaling limit of the particle location  $\{c^{-\beta/\alpha}S_{N_{[ct]}}, t\geq 0\}$ is $\{A_{E_t}, t\geq 0\}$,
  a symmetric stable L\'evy motion time-changed by an inverse stable subordinator.

The random variable $Z_t$ has a smooth density.  For properly scaled waiting times, the density of the standard stable subordinator $Z_t$ has Laplace transform $\E[e^{-\eta Z_t}]=e^{-t\eta^\beta}$ for any $\eta,t>0$, and $Z_t$ is identically distributed with $t^{1/\beta}Z_1$.  Writing $g_\beta(u)$ for the density of $Z_1$, it follows that $Z_s$ has density $s^{-1/\beta}g_\beta(s^{-1/\beta}u)$ for any $s>0$.  Using the inverse relation $ \P(E_t\leq s)= \P(Z_s\geq t)$ and taking derivatives, it follows that $E_t$ has the density
\begin{equation}\begin{split}\label{Etdens}
f_t(s)=\frac d{ds} \P(Z_s\geq t)
&=t\beta^{-1}s^{-1-1/\beta}g_\beta(ts^{-1/\beta}) .
\end{split}\end{equation}
For more details, see \cite{Zsolution,limitCTRW}.

\section{Fractional calculus}

The Caputo fractional derivative of order $0<\beta< 1$, defined by
\begin{equation}\label{CaputoDef}
\frac{\partial^\beta f(t)}{\partial t^\beta}=\frac{1}{\Gamma(1-\beta)}\int_0^t \frac{\partial f(r)}{\partial r}\frac{dr}{(t-r)^\beta},
\end{equation}
was invented to properly handle initial values \cite{Caputo,EIK}.  Its Laplace transform (LT) $s^\beta \tilde f(s)-s^{\beta-1} f(0)$ incorporates the initial value in the same way as the first derivative.  Here $\tilde f(s)=\int_0^\infty e^{-st}f(t)dt$ is the usual Laplace transform.  The Caputo derivative has been widely used to solve ordinary differential equations that involve a fractional time derivative  \cite{GorenfloSurvey,Podlubny}.  In particular, it is well known that the Caputo derivative has a continuous spectrum, with eigenfunctions given in terms of the Mittag-Leffler function
$$E_\beta(z)=\sum_{k=0}^{\infty}\frac{z^k}{\Gamma (1+\beta k)} .$$
In fact, $f(t)=E_\beta(-\lambda t^\beta)$ solves the eigenvalue equation
\[\frac{\partial^\beta f(t)}{\partial t^\beta}=-\lambda f(t)\]
for any $\lambda>0$.  This is easy to check, differentiating term-by-term and using the fact that $t^p$ has Caputo derivative $t^{p-\beta}\Gamma(p+1)/\Gamma(p+1-\beta)$ for $p>0$ and $0<\beta\leq 1$.

For $0<\alpha<2$, the fractional Laplacian  $\Delta^{\alpha/2} f$ is defined for
\[f\in {\rm Dom}(\Delta^{\alpha/2})=
\left\{ f\in L^2(\R^d; dx): \,
\int_{\R^d} |\xi |^\alpha \, | \wh f (\xi)|^2 d\xi <\infty\right\} \]
as the function with Fourier transform
\begin{equation}\label{e:3.2}
  \wh {\Delta^{\alpha/2}f} (\xi)
  =-|\xi|^\alpha \, \wh f(\xi) .
\end{equation}
For suitable test functions (for example, $C^2$ functions with
  bounded second derivatives), the fractional Laplacian can be
  defined pointwise:
\begin{equation}\label{ZQ1}
\Delta^{\alpha/2}f(x)= \int_{y\in\rd} \left( f(x+y)-f(x)
-\nabla f(x)\cdot y {\bf 1}_{\{|y| \leq 1\}}\right)\frac{c_{d,\alpha}}{|y|^{d+\alpha}} \, dy,
\end{equation}
where $c_{d, \alpha}>0$ is a specific constant that depends on $d$ and $\alpha$
so that
$$
c_{d, \alpha} \int_{y\in\rd} \frac{1- \cos y_1}
{|y|^{d+\alpha}} dy=1.
 $$

\begin{remark}\label{R:4.1}
(i) It can be verified using Fourier transforms that,
for $f\in {\rm Dom}(\Delta^{\alpha/2})$, if the right hand side
of \eqref{ZQ1} is well-defined for a.e.\ $x\in \R^d$, then
the Fourier transform of the right-hand side of \eqref{ZQ1}
equals  $-| \xi |^\alpha \wh f(\xi)$
(cf. \cite[Theorem 7.3.16]{RVbook}).
Conversely, it can also be verified that if $f\in L^2(\R^d; dx)$
is a function such that the right hand side
of \eqref{ZQ1} is well-defined for a.e. $x\in \R^d$
and is $L^2(\R^d; dx)$-integrable,
then $f\in {\rm Dom}(\Delta^{\alpha/2})$ and \eqref{ZQ1} holds.

(ii) Using a Taylor series expansion in \eqref{ZQ1},
it is easy to see that
  $\Delta^{\alpha/2} f(x_0)$ exists and is finite at a point
$x_0\in \R^d$ if
$f$ is bounded on $\R^d$ and $f$ is $C^2$ at the point $x_0$.
Hence, if $f$ is bounded and continuous on $\R^d$
 and $f$ is $C^2$ in an open set $D$, then $\Delta^{\alpha/2}f$
 exists pointwise and is continuous in $D$.
 Moreover, if $f$ is a $C^1$ function on $[0, \infty)$ with
 $|f'(t)| \leq c \, t^{\gamma -1}$ for some $\gamma >0$,
 then by \eqref{CaputoDef}, the Caputo fractional derivative
 ${\partial^\beta f (t)}/{\partial t^\beta}$ of $f$ exists
 for every $t>0$ and the derivative is continuous in $t>0$.
 \qed
\end{remark}

For $0<\alpha\leq 2$, let $X$ be the L\'evy process on $\R^d$
such that
 $$ \E\left[ e^{i\xi \cdot (X_t-X_0)}\right]
     =e^{-t|\xi|^\alpha}
 \quad \hbox{for every } \xi \in \R^d.
 $$
 This L\'evy process $X$ is called a standard (rotationally) symmetric $\alpha$-stable process on $\R^d$.
 When $\alpha=2$, it is Brownian motion running at double speed.

Denote the transition semigroup of $X$ by $\{P_t, t>0\}$.
Using the fact that $X_t\Rightarrow X_0$ as $t\to 0+$, it is not hard to show (e.g., see \cite[Theorem 13.4.2]{applebaum}) that $\{P_t, t\geq 0\}$ is a symmetric strongly continuous
 semigroup on the Banach space $L^2(\R^d; dx)$.
Let $(\sF, \sE)$ be the Dirichlet form of $X$ on $L^2(\R^d; dx)$.
That is,
\begin{eqnarray}\label{e:3.4}
\sF &=&
\left\{ u\in L^2(\R^d; dx):
\sup_{t>0} \frac1t (u-P_t u, u)_{L^2(\R^d; dx)} <\infty \right\} , \\
\sE(u, v) &=& \lim_{t\to 0} \frac1t (u-P_t u, v)_{L^2(\R^d; dx)}
\qquad \hbox{for } u, v \in \sF. \label{e:3.5}
\end{eqnarray}
It is known that, for example, via  Fourier transforms \cite{FOT},
\begin{eqnarray*}
\sF &=& W^{\alpha/2, 2}(\R^d):=
 \left\{ u\in L^2(\R^d; dx): \int_{\R^d \times \R^d} \frac{(u(x)-u(y))^2}{|x-y|^{d+\alpha}} dxdy <\infty \right\}, \\
\sE(u, v) &=&  \frac{c_{d, \alpha}}2 \int_{\R^d \times \R^d}
\frac{(u(x)-u(y))(v(x)-v(y))}{|x-y|^{d+\alpha}} dxdy.
\end{eqnarray*}
Let $({\rm Dom}(\L), \L)$ be the $L^2$-generator of the Dirichlet
form $(\sE, \sF)$; that is, $f\in {\rm Dom}(\L) $
if and only if $f\in W^{\alpha/2, 2}(\R^d)$ and
there is some $u\in L^2(\R^d; dx)$ so that
$$ \sE(f, g) = -(u, g) \qquad \hbox{for every }
g\in W^{\alpha/2, 2}(\R^d);
$$
  in this case, we denote this $u$ by $\L f$. It is known
(cf. \cite{FOT})
that $\L$ is also the semigroup generator of $\{P_t, t>0\}$
on the space $L^2(\R^d; dx)$. Using the Fourier transform,
one can conclude (cf. \cite{FOT}) that
$f\in {\rm Dom}(\L)$ if and only if
$\int_{\R^d} |\xi |^\alpha | \wh f (\xi)|^2 d\xi <\infty$,
and $\wh {\L f} (\xi)=-|\xi|^\alpha \wh f (\xi)$
for every $f\in {\rm Dom} (\L)$.  Hence the $L^2$-generator of $X$
is the fractional Laplacian $\Delta^{\alpha/2}$.

It follows directly from Dirichlet form theory (cf.
\cite{FOT}) that,  for $f\in L^2(\R^d)$ and $t>0$,
$P_t f \in \sF=W^{\alpha/2, 2}(\R^d)$,
and  $v(t, x):=\E_x[f(X_t)]$ is a weak solution to the following
  parabolic equation:
\begin{equation}\label{spaceFDE}
\frac{\partial}{\partial t} v(t, x)=\Delta^{\alpha/2}  v(t, x);
\quad v(0, x)=f(x) .
\end{equation}
That is, the function $x\mapsto v(x,t)$ belongs to the domain of the $L^2$ generator $\L=\Delta^{\alpha/2}$ for every $t>0$, and equation \eqref{spaceFDE} holds in the space $L^2(\R^d;dx)$.  Here the fractional Laplacian and the first time derivative in \eqref{spaceFDE} are defined in terms of the Banach space norm.  For example, the time derivative is the limit of a difference quotient that converges in the $L^2$ sense, so it need not exist point-wise.   The classical diffusion equation models the evolution of particles away from their starting point, due to molecular collisions.  The space-fractional diffusion equation \eqref{spaceFDE} models particle motions in a heterogeneous environment, where the probability of long particle jumps follows a power law \cite{fadePRE}.

For $0<\alpha <2$, the symmetric $\alpha$-stable process $X$
  can be obtained from Brownian motion
 on $\R^d$ through subordination in the sense of Bochner \cite{Bochner}.
 Let $\{B , \P_x, x\in \R^d\}$ be Brownian motion on $\R^d$
 with $\P_x(B_0=x)=1$ and  $\E_0 [B_t B_t ']=2tI$, where $'$ denotes the transpose, and $I$ is the $d\times d$ identity matrix.
For $0<\alpha<2$, let $Z_t$ be a standard stable subordinator with $Z_0=0$,
whose Laplace transform is $\E [ e^{-s Z_t}]=e^{-t s^{\alpha/2}}$
for every $s, t>0$. Then it is easy to verify, using Fourier transforms and a simple conditioning argument, that $B_{Z_t}$ is a symmetric $\alpha$-stable L\'evy process
starting from the origin
that has the same distribution as $X$, with $X_0=0$.
The process $X$  has
a jointly continuous  transition density function $p(t,x,y)=p_t(x-y)$ with respect to the Lebesgue measure
in $\R^d$.
That is,
$$
\P_x (X_t\in A )=\int_A p(t,x,y)dy.
$$
Using the self-similarity of the stable process and
its relation with Brownian motion through subordination, it is not hard to show that for $\alpha\in(0,2)$ we have
\begin{equation}\label{ultraconractive-bound}
p_t(x)=t^{-d/\alpha}p_1(t^{-1/\alpha}x)\leq t^{-d/\alpha}p_1(0)=:t^{-d/\alpha}M_{d,\alpha},
\quad t>0, x\in \rd.
\end{equation}

Another kind of time change relates to particle waiting times.  Suppose $\{T_t, t\geq 0\}$ is a uniformly bounded strongly continuous semigroup  on a Banach space $E$, with infinitesimal generator $(\A, \operatorname{Dom}(\A))$.
It is known that $v(t)=T_tf$ solves the Cauchy problem $\partial v/\partial t=\A v$ with $v(0)=f$ for any $f\in \operatorname{Dom}(\A)$
(see \cite{ABHN}). Let $Z$ be a standard $\beta$-stable subordinator independent of $X$, and recall that $E_t=\inf\{s>0: Z_s>t\}$ is its inverse process.  If $g_\beta (u)$ is the density of $Z_1$, then \cite[Theorem 3.1]{fracCauchy} shows that another subordinated semigroup
\begin{equation}\label{e:3.8}
R_tf=\int_0^\infty  g_\beta(u) T_{(t/u)^\beta}f\,du
\end{equation}
yields solutions to the time-fractional Cauchy problem:  $w(t)=R_tf$ solves
\[\frac{\partial^\beta}{\partial t^\beta}w(t)=\A w;\quad w(0)=f\]
on the Banach space $E$ for any $f\in \operatorname{Dom}(\A)$. Applying this to the transition semigroup $\{P_t, t\geq 0\}$ of the symmetric $\alpha$-stable process $X$ on the space $L^2(\rd;dx)$, one sees that the process $Y_t=X_{E_t}$ can be used to solve the space-time diffusion equation
on $\R^d$; that is,
$w(t, x)=\E_x[f(Y_t)]$ is a weak solution for
\begin{equation}\label{space-timeFDE}
\frac{\partial^\beta}{\partial t^\beta} w(x,t)= \Delta^{\alpha/2}  w(x,t);\quad w(x,0)=f(x) .
\end{equation}
That is, the function $x\mapsto w(x,t)$ belongs to the domain of the $L^2$ generator $\L=\Delta^{\alpha/2}$ for every $t>0$, and equation \eqref{space-timeFDE} holds in the Banach space $L^2(\rd;dx)$.

\section{Eigenfunction expansion for bounded domains}

Let $D $ be a  bounded open subset of $\RR{R}^d$.
Recall that $X$ is a standard spherically symmetric stable process
 on $\rr^d$, and define the first exit time
\[\tau_D =\inf \{ t\geq 0:\ X_t\notin D\} .\]
Let $X^D$ denote the process $X$ killed upon leaving $D$;
that is, $X^D_t=X_t$ for $t<\tau_D$ and $X^D_t=\partial$ for
$t\geq \tau_D$. Here $\partial$ is a cemetery point added to $D$.
Throughout this paper, we use the convention that any real-valued
function $f$ can be extended by taking
$f(\partial )=0$.
The subprocess $X^D$ has a jointly continuous transition density
function $p_D(t, x, y)$ with respect to the Lebesgue measure
on $D$. In fact, by the strong Markov property of $X$, one has
for $t>0$ and $x, y\in D$,
\begin{equation}\label{e:4.1a}
 p_D(t, x, y)=p(t, x, y)-\E_x [ p(t-\tau_D, X_{\tau_D}, y);
\tau_D<t] \leq p(t, x, y).
\end{equation}
Denote by $\{P^D_t, t\geq 0\}$ the transition semigroup of
$X^D$, that is
$$
P^D_t f(x)=\E_x[f(X^D_t)]= \int_D p_D(t, x, y) f(y) dy.
$$
The proof of the following facts can be found in \cite{FOT}:
The operators $\{P^D_t, t\geq 0\}$ form a
symmetric strongly continuous contraction semigroup in $L^2(D; dx)$.
Let $(\sE^D, \sF^D)$ denote the Dirichlet form of $X^D$,
defined by \eqref{e:3.4}--\eqref{e:3.5} but with
 $\{P^D_t, t>0\}$ in place of $\{P_t, t>0\}$.
 Then $\sF^D$ is the $\sqrt{\sE_1}$-completion of the space
$C^\infty_c(D)$ of smooth functions with compact support in $D$,
denoted by $W^{1,2}_0(D)$ in literature. Here
$\sE_1 (u, u)=\sE (u, u)+\int_{\R^d} u(x)^2 dx$. Moreover,
$\sE^D(u, v)=\sE(u, v)$ for $u, v\in W^{\alpha/2, 2}_0(D)$.
Let $\L_D$ be the $L^2$-infinitesimal generator of
$(\sE^D, \sF^D)$; that is, its
domain ${\rm Dom}(\L_D)$ consists
all $f\in W^{\alpha/2, 2}_0(D)$ such that
$$ \sE^D(f, g)=-(u, g)_{L^2(D; dx)} \qquad \hbox{for every }
g\in W^{\alpha/2, 2}_0(D);
$$
for some $u\in L^2(D; dx)$; in this case, we denote this $u$ by $\L_D f$.
It is well-known (cf. \cite{FOT}) that $\L_D$ is the $L^2$-generator
of the strongly continuous semigroup $\{P^D_t, t>0\}$ in $L^2(D; dx)$.
For every $f\in L^2(D; dx)$ and $t>0$,
$P^D_t f \in {\rm Dom}(\L_D)\subset W^{\alpha/2, 2}_0(D)$.
Moreover $u(t, x):=P^D_t f (x)$
is the unique weak solution to
$$ \frac{\partial u}{\partial t} = \L_D u
$$
with initial condition $u(0, x)=f(x)$ on the Banach space $L^2(D; dx)$.

Note that the transition kernel $p_D(t,x,y)$ is symmetric
and strictly positive   with
\begin{equation}\label{star}
p_D(t,x,y)\leq p(t,x,y)\leq t^{-d/\alpha}M_{d,\alpha}, \ \ x,y\in D, t>0
\end{equation}
in view of \eqref{ultraconractive-bound}.
In particular, one has $\sup_{x\in D} \int_D p(t,x, y)^2 dy<\infty$
for every $t>0$.
Thus for each $t>0$,
$P^D_t$ is a Hilbert-Schmidt operator in $L^2(D; dx)$
so it is compact.
Therefore there is a sequence of positive numbers
 $0<\lambda_1<\lambda_2\leq \lambda_3\leq\cdots$
 and an orthonormal basis
$\{\psi_n, n\geq 1\}$ of  $L^2(D; dx)$ so that
$ P^D_t \psi_n =e^{-\lambda_n t} \psi_n$ in $L^2(D; dx)$
for every $n\geq 1$ and $t>0$.
Since for every $f\in L^2(D; dx)$,
$f(x)= \sum_{n=1}^\infty \< f, \psi_n\> \psi_n (x)$, we have
\begin{equation}\label{eigen-f}
 P^D_t f (x) =\sum_{n=1}^\infty \< f, \psi_n\> P^D_t \psi_n (x)
= \sum_{n=1}^\infty e^{-\lambda_n t } \< f, \psi_n\>   \psi_n (x).
\end{equation}
That is, the transition density
\begin{equation}\label{eigen-kernel}
p_D(t,x,y)=\sum_{n=1}^\infty e^{-\lambda_n t}\psi_n(x)\psi_n(y) .
\end{equation}
It follows from \cite[Theorem 2.3]{BG} that for any bounded open
subset $D$ of $\R^d$, one has
\begin{equation}\label{e:4.2}
c_1 n^{\alpha/d} \leq \lambda_n \leq c_2 n^{\alpha/d}
\qquad \hbox{for every } n\geq 1.
\end{equation}
Using the spectral representation, one has
\begin{equation}\label{domain-def}
\operatorname{Dom} (\L_D)=
\left\{f\in L^2(D): \ \ \| \L_D f \|_{L^2(D)}^2=\sum_{n=1}^\infty \lambda_n^2\l f, \psi_n \r ^2<\infty \right\}.
\end{equation}
and
$$ \L_D f (x) = - \sum_{n=1}^\infty \lambda_n \< f, \psi_n\>
\psi_n(x) \qquad \hbox{for } f\in {\rm Dom} (\L_D).
$$
For any real valued function $\phi: \R \to \R$, one can also define
the operator $\phi (\L_D)$ as follows:
\begin{eqnarray*}
{\rm Dom}(\phi (\L_D))&=& \left\{f\in L^2(D; dx): \sum_{n=1}^\infty
\phi (\lambda_n)^2 \<f, \psi_n\>^2 <\infty \right\}, \\
\phi(\L_D) f &=& \sum_{n=1}^\infty \phi (\lambda_n) \< f , \psi_n\> \psi_n.
\end{eqnarray*}
In next section, the operator $\L_D^k$ defined using $\phi(t)=t^k$ will be utilized.

\medskip
The generator $\L_D$
is also called the fractional Laplacian on $D$ with zero exterior
condition, denoted as $\Delta^{\alpha/2}|_D$.
We  now record a lemma that gives an explicit expression of
$\L_D$.

\begin{lemma}\label{L:4.1}
 For $f\in\sF^D$, if
\begin{equation}\label{e:4.8}
 \phi (x):= \lim_{\eps \to 0} \int_{\{y\in \R^d: |y-x|>\eps \}}
(f(y)-f(x)) \frac{c_{d, \alpha}}{|y-x|^{d+\alpha}} dy
\end{equation}
 exists and the convergence is uniformly on each compact subsets of $D$
and $\phi \in L^2(D; dx)$, then $f\in {\rm Dom}(\L_D)$
and $\phi = \L_D f$. In particular, if $f$ is a bounded function in $\sF^D\cap C^2(D)$, then $f\in {\rm Dom}(\L_D)$ and
\begin{eqnarray*}
\L_D f(x) &=&  \lim_{\eps \to 0} \int_{\{y\in \R^d: |y-x|>\eps \}}
(f(y)-f(x)) \frac{c_{d, \alpha}}{|y-x|^{d+\alpha}} dy \\
&=& \int_{y\in \R^d} \left(f(x+y)-f(x) -\nabla f(x) \cdot y {\bf 1}_{\{|y|\leq 1\}}\right) \frac{c_{d, \alpha}}{| y|^{d+\alpha}} dy.
\end{eqnarray*}
\end{lemma}

\proof Suppose that $f\in \sF^D$ and that $\phi$ defined by
\eqref{e:4.8} converges locally uniformly in $D$ and is
in  $L^2(D; dx)$. Then for every $g\in C^2_c(D)$, by the
expression of $\sE^D(f, g)$ and the symmetry,
\begin{eqnarray*}
\sE^D(f, g) &=& \frac12  \int_{  \R^d\times \R^d }
(f(x)-f(y))(g(x)-g(y))\frac{c_{d, \alpha}}{|x-y|^{d+\alpha}} dx dy \\
 &=& \frac12 \lim_{\eps \to 0}
\int_{\{(x, y)\in \R^d\times \R^d: |x-y|>\eps\}}
(f(y)-f(x))(g(y)-g(x))\frac{c_{d, \alpha}}{|y-x|^{d+\alpha}} dx dy \\
&=& - \lim_{\eps \to 0}
\int_{\R^d} \left(\int_{\{y\in \R^d: |y-x|>\eps\}}
(f(y)-f(x))\frac{c_{d, \alpha}}{|y-x|^{d+\alpha}} dy \right)
g(x) dx \\
&=& - \int_{\R^d} \phi (x) g(x) dx.
\end{eqnarray*}
Since $C^2_c(D)$ is $\sE^D_1$-dense in $W^{\alpha/2, 2}_0(D)$,
this implies that $f\in {\rm Dom}(\L_D)$ and $\L_D f =\phi$
on $D$.

Assume now that  $f$ is a bounded function in $\sF^D\cap C^2(D)$.
Using a Taylor expansion, one easily sees that
$$\int_{y\in \R^d} \left| f(x+y)-f(x) -\nabla f(x) \cdot y {\bf 1}_{\{|y|\leq 1\}}\right| \frac{c_{d, \alpha}}{| y|^{d+\alpha}} dy
<\infty \qquad \hbox{for every } x\in D
$$
and the integral is a continuous function on $D$.
Set
$$
\psi (x)= \int_{y\in \R^d} \left(f(x+y)-f(x) -\nabla f(x) \cdot y {\bf 1}_{\{|y|\leq 1\}}\right) \frac{c_{d, \alpha}}{| y|^{d+\alpha}} dy
\qquad \hbox{for } x\in D.
$$
For any  compact subset $K$ of $D$, let
$$
K_\eps:=\{z\in \R^d:
\hbox{there is some } x\in K \hbox{ so that } |z-x|\leq \eps\}.
$$
Defining
\[\| D^2 f\|_\infty=\max_{1\leq i, j\leq d} \left\| \frac{\partial^2 f}{\partial x_i
\partial x_j}\right\|_\infty, \]
 we have
\begin{eqnarray*}
&& \lim_{\eps\to 0} \sup_{x\in K} \left|
\int_{\{y\in \R^d: |y-x|>\eps \}}
(f(y)-f(x)) \frac{c_{d, \alpha}}{|y-x|^{d+\alpha}} dy -\psi (x)
\right| \\
&=& \lim_{\eps\to 0} \sup_{x\in K} \left|
\int_{\{y\in \R^d: |y-x|\leq \eps \}}
  \left(f(x+y)-f(x) -\nabla f(x) \cdot y {\bf 1}_{\{|y|\leq 1\}}\right) \frac{c_{d, \alpha}}{| y|^{d+\alpha}} dy
\right| \\
&\leq &  \lim_{\eps\to 0}   \left|
\int_{\{y\in \R^d: |y-x|\leq \eps \}}
  \sup_{z\in K_\eps} \| D^2 f\|_\infty \, |y|^2 \,
  \frac{c_{d, \alpha} }{| y|^{d+\alpha}} dy
\right| =0 .
\end{eqnarray*}
By what we have shown in the first part, this implies
 that $f\in {\rm Dom}(\L_D)$ with $\L_D f=\psi$,
  which completes the proof of the lemma.
\qed

\bigskip

The main purpose of this paper is to investigate the existence
of strong solution to the following equation:
\begin{equation}\begin{split}\label{e:4.3}
\frac{\partial^\beta}{\partial t^\beta}u(t,x) &=
\Delta^{\alpha/2} u(t,x); \quad  x\in D, \ t>0\\
u(t,x)&=0, \quad  x\in  D^c, \ t>0, \\
u(0,x)& = f(x), \quad x\in D.
\end{split}\end{equation}
Let $C_\infty(D)$ denote the Banach space of bounded continuous functions on $\rr^d$ that vanish off $D$, with the sup norm.

\begin{defi} (i) Suppose that $f\in L^2(D; dx)$. A function $u(t, x)$ is said to be a weak solution to
\eqref{e:4.3} if  $u(t, \cdot) \in W^{1,2}_0(D)$ for every $t>0$,
  $\lim_{t\downarrow 0} u(t, x) =f(x)$ a.e.\ in $D$,
  and  ${\partial^\beta}/{\partial t^\beta}u(t,x) =
\Delta^{\alpha/2} u(t,x)$
in the distributional sense; that is,
for every $\psi \in C^1_c([0, \infty)$ and $\phi \in C^2_c(D)$,
$$ \int_{\R^d} \left( \int_0^\infty u(t, x) \frac{\partial^\beta \psi (t)}{\partial t^\beta} dt \right) \phi (x)  dx
= \int_0^\infty \sE^D( u (t, \cdot ), \phi ) \, \psi (t) \, dt .
$$

(ii) Suppose that $f\in C(D)$.
 A function $u(t, x)$ is said to be a strong solution
\eqref{e:4.3} if for every $t>0$, $u(t, \cdot)\in C_\infty (D)$,
$\Delta^{\alpha/2} u(t,\cdot )(x)$ exists pointwise for every
$x\in D$ in the sense of \eqref{ZQ1}, the Caputo fractional
derivative ${\partial^\beta u(t, x)}/{\partial t^\beta}$
exists pointwise for every $t>0$ and $x\in D$,
${\partial^\beta}/{\partial t^\beta}u(t,x) =
\Delta^{\alpha/2} u(t,x)$ pointwise in $(0, \infty)\times D$,
and $\lim_{t\downarrow 0} u(t, x) =f(x)$ for every $x\in D$.
 \end{defi}

\medskip

A boundary point $x$ of an open set  $D$ is said to be {\it regular}
 for $D$ if ${\mathbb P}_x[\tau_D(X)=0]=1$.
 A sufficient condition for $x_0\in \partial D$ to be  regular
  for $D$ is
  that $D$ satisfies an {\em exterior cone condition} at $x_0$, that is,   there exists a finite right
circular open cone $V=V_{x_0}$ with vertex $x_0$ such that
$V_{x_0} \subset D^c$ (cf. \cite[Theorem 2.2]{ChenSong97}). An open set $D$ is said to be regular
if every boundary point of $D$ is regular for $D$.
Assume now that $D$ is a regular open set.
Then \cite[Theorem 2.3]{ChenSong97}
 shows that $\{P^D_t, t>0\}$ is a strongly continuous (Feller) semigroup on the Banach space $C_\infty(D)$ of bounded continuous functions on $\rr^d$ that vanish off $D$, with the sup norm.  Moreover, $\{P^D_t, t>0\}$ has the same set of  eigenvalues and eigenfunctions on $C_\infty(D)$ as on $L^2(D; dx)$: $P^D_t\psi_n=e^{-\lambda_n t}\psi_n$ in $C_\infty(D)$ (see \cite[Theorem 3.3]{ChenSong97}).     In particular, every eigenfunction $\psi_n$ of the $L^2$-generator $\L_D$ is a bounded continuous function on $D$ that vanishes continuously on
 the boundary $\partial D$.

\medskip

\section{Space-time fractional diffusion in bounded domains}

In this section, we
prove  strong  solutions to space-time fractional diffusion equations
on bounded domains in $\rd$.  We give an explicit solution
formula, based on the solution of the corresponding Cauchy
problem.  The basic argument uses an eigenfunction
expansion of the fractional Laplacian on $D$, and separation of variables.  The probabilistic representation of the solution is constructed from a killed stable processes, whose index corresponds to the fractional Laplacian, modified by an inverse stable time change, whose index equals the order of the fractional time derivative.

Recall that $X$ is a rotationally symmetric $\alpha$-stable process
in $\R^d$ and   $\{E_t, t\geq 0\}$ is the inverse of a standard stable subordinator of index $\beta\in (0,1)$, independent of $X$.  In the following proof, we denote by $c,c_1,c_2,\ldots$ a constant that may change from line to line.

\begin{theorem}\label{thm:5.6}
Let $D$ be a regular open subset of $\R^d$.
Suppose  $f\in {\rm Dom}(\L_D^k)$ for some
 $k>-1+({3d+4})/({2\alpha})$.
Then \[u(t, x)=\E_x [ f(X^D_{E_t})]\in C_b([0, \infty) \times \R^d)
\cap C^{1,2}((0, \infty)\times D)\]
and $u(t,x)$ is a strong solution to the space-time fractional diffusion equation \eqref{e:4.3}.
\end{theorem}

\proof First we will prove that  $f\in C_\infty (D)$.
Let $0<\lambda_1<\lambda_2 \leq \lambda_3 \leq \cdots$
be the eigenvalues of $\L_D$ and $\{\psi_n, n\geq 1\}$
be the corresponding eigenfunctions, which form
an orthonormal basis for $L^2(D; dx)$.
Note that, since $D$ is a regular open set,
we have from the last section that  $\psi_n\in C_\infty (D)$
for each $n\geq 1$. Since $f\in {\rm Dom}(\L_D^k)$ for some
 $k>-1+({3d+4})/({2\alpha})$, using \eqref{e:4.2} it follows that
\begin{equation}\label{e:5.6}
 M:=\sum_{n=1}^\infty \lambda_n^{2k}\<f, \psi_n\>^2 <\infty,
\end{equation}
and so $|\<f, \psi_n\>| \leq \sqrt{M} \lambda_n^{-k}$.  From \eqref{star}
and \eqref{eigen-kernel}  we get
$$e^{-\lambda_n t}|\psi_n(x)|^2\leq \sum_{k=1}^\infty e^{-\lambda_k t} |\psi_k (x)|^2=p_D(t,x,x)\leq M_{d,\alpha} t^{-d/\alpha}$$
and hence, taking square roots of both sides, we get
$$
|\psi_n(x)|\leq e^{\lambda_n t/2} \sqrt{M_{d,\alpha} t^{-d/\alpha}}
$$
Taking $t=1/\lambda_n$ gives us
\begin{equation}\label{e:5.2}
|\psi_n (x)| \leq c \lambda_n^{d/(2\alpha)}
\qquad \hbox{for every } x\in D
\end{equation}
for some $c>0$.  Since $k>-1+({3d+4})/({2\alpha})$, \eqref{e:5.2} together with
 \eqref{e:4.2} implies that
$$\sum_{n=1}^\infty | \< f, \psi_n\>| \, \| \psi_n\|_\infty
\leq c \sum_{n=1}^\infty \lambda_n^{-k} \lambda_n^{d/(2\alpha)}
\leq c \sum_{n=1}^\infty n^{(\alpha/d)(d/(2\alpha)-k)}<\infty.
$$
Hence $f(x)=\sum_{n=1} \<f, \psi_n\> \psi_n$
converges uniformly on $D$, and so $f\in C_\infty (D)$.

Recall that $P^D_tf(x)=\E_x [ f(X^D_t) ]$
is the unique weak solution in $W^{\alpha/2, 2}_0(D)$
of the equation
\begin{equation}\label{e:5.27}
\frac{\partial}{\partial t} v(t, x)=\Delta^{\alpha/2}v(t, x)\ \quad \hbox{ with }  v(0, x)=f(x)
\end{equation}
on the Banach space $L^2(\rd;dx)$ (cf. (see \cite{FOT}).
The semigroup $P^D_t$ has density function $p_D(t, x, y)$ given by
\eqref{e:4.1a}.
Note that  $p(t, x, y)$ is smooth in $x$. By a proof similar to  \cite[Proposition 3.3]{BC}, we have for every $j\geq 1$ and
$1\leq i\leq d$ that
\begin{equation}\label{e:5.29}
 \Big| \frac{\partial^j}{\partial x_i^j} p(t, x, y) \Big|
  \leq c\left( t^{-(d+j)/\alpha} \wedge \frac{t} {|x-y|^{d+\alpha+j}}\right)
  \leq c_1 t^{-j/\alpha} p(t, x, y).
\end{equation}
In view of the symmetry $p(t, x, y)=p(t, y, x)$ and $p_D(t, x, y)
=p_D(t, y, x)$, we have
  from \eqref{e:4.1a} and \eqref{e:5.29}
that $P^D_tf(x)=\int_D p_D(t, x, y)f(y)dy$ is smooth in $x\in D$. Moreover,
for every compact subset $K$ of $D$ and $T>0$, there is a constant $c_2=c_2(d, \alpha, K, T)$ such that, for $x\in K$ and $t\in (0,  T]$,
\begin{equation}\label{e:5.4}
 \Big| \frac{\partial^j}{\partial x_i^j} p_D(t, x, y) \Big|
  \leq c_2 t^{-j/\alpha} p(t, x, y).
\end{equation}
The Chapman-Kolmogorov equation implies
$$
\int_{\R^d} p(t,x,y)^2dy=\int_{\R^d} p(t,x,y)p(t,y,x)\,dy=p(2t,x,x).
$$
It then follows using \eqref{star}, \eqref{e:5.4}, and the
Cauchy-Schwarz inequality that
\begin{equation}\begin{split} \label{e:5.30}
 |\nabla^jP^D_t f(x)|
\leq  c_3 t^{-j/\alpha} (2t)^{-d/(2\alpha)}\| f\|_{L^2(D)} .
\end{split}\end{equation}
Consequently, each eigenfunction $\psi_n(x)=e^{\lambda_n t} P^D_t \psi_n(x)$ is smooth inside $D$ with
$$ \Big| \nabla^j \psi_n(x)\Big| \leq
c_3 t^{-(d+2j)/(2\alpha)} e^{\lambda_n t}
$$
for $x\in K$ and $t\in (0, T]$. Taking $t=1/\lambda_n$ yields
\begin{equation}\label{e:5.31}
 \Big| \nabla^j \psi_n(x)\Big| \leq
c_3 \lambda_n^{(d+2j)/(2\alpha)} \qquad \hbox{for } x\in K.
\end{equation}
In view of \eqref{eigen-f}, $P^D_tf(x)$ is also differentiable in $t>0$.
(The eigenfunction expansion \eqref{eigen-f} together with \eqref{e:5.31}
 gives another proof that $P^D_tf$ is $C^\infty$ in $x\in D$.)
Hence in view of Remark \ref{R:4.1},
$v(t, x)=P^D_tf(x)$ is a classical
solution for ${\partial v}/{\partial t} = \L_D v$ in $D$.

Now define
$$u(t, x)=\E_x[ f(X^D_{E_t})]  = \E_x[f(X_{E_t});  E_t<\tau_D ].
$$
Since $X^D$ generates a strongly continuous (Feller) semigroup on $C_\infty(D)$,
$P^D_t f(x)$ is a bounded continuous
function on $[0, \infty) \times \R^d$ that vanishes on $[0, \infty)
\times D^c$, and hence so is $u$, in view of  \eqref{e:3.8}.
 By \cite[Theorem 3.1]{fracCauchy}
(and  \cite[Theorem 4.2]{limitCTRW}), $u(t, x)$ is a weak solution
for the parabolic equation \eqref{e:4.3} on $L^2(\rd,dx)$.
Then, to show that $u$ is a classical solution, by Remark \ref{R:4.1},
it suffices to show that
$u(t, \cdot)\in C^2( D)$ for each $t>0$, and that the Caputo derivative
of $t\mapsto u(t, x)$ exists for each $x$, and is jointly continuous in $(t, x)$.

Bingham \cite{Bingham} showed that the inverse stable law $E_t$ with density $f_t(s)$ given by \eqref{Etdens} has a Mittag-Leffler distribution, with Laplace transform $\E[e^{-\lambda E_t}]=E_\beta(-\lambda t^\beta).$
Then it follows, using \eqref{eigen-f} and a simple conditioning argument, that
\begin{equation}\label{mittag-series-2}
\begin{split}
u(t,x)
&=\int_0^{\infty}\E_x \left[ f(X_s); s<\tau_D \right] f_t(s)\, ds\\
&=\int_0^{\infty}\left( \sum_{n=1}^\infty e^{-s\lambda_n} \langle f,\psi_n\rangle \psi_n(x)
 \right) f_t(s)\, du\\
&=\sum_{n=1}^\infty E_\beta(-\lambda_n t^\beta) \langle f,\psi_n\rangle \psi_n(x) .
\end{split}
\end{equation}
Then, since $0\leq E_\beta (-\lambda_n t^\beta)\leq {c}/({1+\lambda_n t^\beta})$, we have by \eqref{e:5.31} and \eqref{mittag-series-2} that
\begin{eqnarray*}
\|\nabla^j u\|_\infty
&\leq &\sum_{n=1}^\infty E_\beta(-\lambda_nt^\beta) |\<f, \psi_n\>|\,
\|\nabla^j \psi_n\|_\infty\\
&\leq& \sum_{n=1}^\infty c    \, \lambda_n^{-k}\sqrt{M}\,
\frac{\lambda_n^{(d+4)/(2\alpha)}}{1+\lambda_n t^\beta} \\
&\leq & (c\sqrt{M}) t^{-\beta} \sum_{n=1}^\infty  \lambda_n^{(d+4)/(2\alpha)-1-k}
\end{eqnarray*}
 for $j=1, 2$.  Then by \eqref{e:4.2},
\begin{eqnarray*}
\|\nabla^j u\|_\infty
&\leq & (c\sqrt{M}) t^{-\beta} \sum_{n=1}^\infty  \lambda_n^{(d+4)/(2\alpha)-1-k}\\
&\leq &(c c_\alpha \sqrt{M}) t^{-\beta}\sum_{n=1}^\infty  n^{({\alpha}/{d})((d+4)/(2\alpha)-1-k)}<\infty
\end{eqnarray*}
if $k> ({3d+4-2\alpha})/({2\alpha})$.
This proves that, when $k>-1+({3d+4})/({2\alpha})$,
$u(t, x)$ is $C^2$ in $x\in K$, and hence in $D$.
Consequently, by Remark \ref{R:4.1}, the spatial
fractional derivative $\Delta^{\alpha/2} u(t, x)$
exists pointwise for $x\in D$,  and
is a jointly continuous function in $(t, x)$.

Next we show $u(t, x)$ is $C^1$ in $t>0$.
Let $0<\gamma <1\wedge ({4}/({2\alpha}) -1)$.
By \cite[Equation (17)]{krageloh},
$$
\left| \frac{\partial }{\partial t}E_\beta(-\lambda_n t^\beta)\right|\leq c\frac{\lambda_n t^{\beta-1}}{1+\lambda_n t^\beta}
\leq c \, \lambda_n^\gamma \, t^{\gamma \beta -1}.
$$
This together with  \eqref{e:5.6} and \eqref{e:5.2}
 yields that
\begin{equation*}\begin{split}
 \sum_{n=1}^\infty \Big| \frac{\partial }{\partial t}E_\beta(-\lambda_n t^\beta) \, \< f, \psi_n\> \psi_n(x) \Big|
 &\leq  \sum_{n=1}^\infty c\lambda_n^\gamma t^{\beta -1} \lambda_n^{-k} \lambda_n^{d/(2\alpha)} \\
&\leq c t^{\gamma \beta -1} \sum_{n=1}^\infty
n^{(\alpha/d ) (\gamma -k +  d/(2\alpha))}
\leq c \,  t^{\gamma \beta -1} .
\end{split}\end{equation*}
Then it follows by a dominated convergence argument
that $u(t, x)$ is continuously differentiable in $t>0$, with
\begin{equation}
\Big| \frac{\partial u(t, x)}{\partial t}\Big|
\leq \sum_{n=1}^\infty  \Big| \frac{\partial }{\partial t}E_\beta(-\lambda_n t^\beta) \, \< f, \psi_n\> \psi_n(x) \Big|
<c t^{\gamma \beta -1} \qquad \hbox{for every } x\in D.
\end{equation}
Hence by Remark \ref{R:4.1},
The Caputo fractional derivative ${\partial^\beta  u(t, x)}/{\partial t^\beta}$ of $u(t, x)$ exists pointwise and is jointly continuous
in $(t, x)$.
Since $u(t, x)$ is a weak solution of \eqref{e:4.3} on $L^2(\rd;dx)$,
by the above regularity property of $u(t, x)$, it is also a strong solution of \eqref{e:4.3}.
\qed

\bigskip

\begin{remark}
The above proof can be easily modified to show that,
if  $D$ is a bounded regular open subset of $\R^d$
  and $f\in {\rm Dom}(\L_D^k)$ for some
 $k>1+ ({3d})/({2\alpha})$, then
 $ u(t, x)=\E_x [ f(X^D_{E_t})] $ is
a weak solution to the space-time fractional diffusion equation \eqref{e:4.3}. Moreover, the Caputo derivative ${\partial^\beta u}/{\partial t^\beta}$ exists
pointwise as a jointly continuous function in $(t, x)$, and
$\L_D u$ has a continuous version that equals ${\partial^\beta u}/{\partial t^\beta}$ on $(0, \infty)\times D$.
\end{remark}
\bigskip

\begin{remark}
The paper \cite{DOFCP} solves distributed-order time-fractional diffusion equations $\partial_t^\nu u=\Delta u$ on bounded domains.  The distributed-order time-fractional derivative is defined by
\[\partial_t^\nu f(t)=\int \frac{\partial^\beta f(t)}{\partial t^\beta} \nu(d\beta) ,\]
where $\nu$ is a positive measure on $(0,1)$.  It may also be possible to extend the results of this paper to develop strong solutions and probabilistic solutions for $\partial_t^\nu u=\Delta^{\alpha/2} u$ on bounded domains.  Distributed-order time-fractional diffusion equations can be used to model ultraslow diffusion, in which a cloud of particles spreads at a logarithmic rate, also called Sinai diffusion \cite{M-S-ultra}. \qed
\end{remark}

\begin{remark}
The fractional Laplacian generates the simplest non-Gaussian stable process in $\R^d$.  Stable processes are useful in applications because they represent universal random walk limits.  For random walks with strongly asymmetric jumps, a wide variety of alternative limit processes exists, see for example \cite{RVbook}.  Because the generators of these processes are not self-adjoint, the extension of results in this paper to that case remains a challenging open problem.
\qed \end{remark}


\begin{thebibliography}{99}

\bibitem{applebaum} Applebaum, D. (2004). \emph{L\'{e}vy Processes and Stochastic Calculus.} Cambridge studies in advanced mathematics.

\bibitem{ABHN}Arendt, W., Batty, C., Hieber, M. and Neubrander, F. (2001). {\it  Vector-valued Laplace transforms and Cauchy problems.} Monographs in Mathematics, Birkh\"{a}user-Verlag, Berlin.

\bibitem{fracCauchy} Baeumer, B. and Meerschaert, M.M. (2001).
Stochastic solutions for fractional Cauchy problems,
\emph{Fractional Calculus Appl. Anal.} {\bf 4} 481--500.


\bibitem{bass}Bass, R. F. (1998). {\it Diffusions and Elliptic Operators.} Springer-Verlag, New York.

\bibitem{BC}  Bass, R.F. and  Chen, Z.-Q. (2006),
Systems of equations driven by stable processes.
{\it Probab. Theory Relat. Fields, \bf 134}, 175-214.

\bibitem{Bingham}Bingham,  N.H. (1971). Limit theorems for occupation times of Markov processes.
  {\it Z. Warsch. verw. Geb.} {\bf 17}, 1--22.

\bibitem{BG} Blumenthal, R.M.  and  Getoor, R.K. (1959).
Asymptotic distribution of the eigenvalues for a class of
Markov operators. {\it Pacific J. Math. \bf 9},
399--408.

\bibitem{Bochner} Bochner, S. (1949). Diffusion equations and stochastic processes. {\it Proc. Nat. Acad. Sci. USA} {\bf 85},
369--370.

 \bibitem{Caputo}  Caputo, M. (1967). Linear models of dissipation whose Q is almost frequency independent, Part II. {\it Geophys. J. R. Astr. Soc.} {\bf 13} 529-539.

\bibitem{ChenSong97}    Chen, Z.-Q. and Song, R.  (1997). Intrinsic ultracontractivity and conditional gauge for symmetric
stable processes, {\it J. Funct. Anal.} {\bf 150}, 204--239.

\bibitem{davies} Davies, E.B. (1989). {\it Heat Kernels and Spectral Theory.} Cambridge Univ. Press, Cambridge.

\bibitem{davies-07} Davies, E.B. (2007). {\it Linear Operators and Their Spectra}. Cambridge Univ. Press, Cambridge.

\bibitem{EIK}Eidelman, S.D., Ivasyshen, S.D., and Kochubei, A.N. (2004). {\it Analytic Methods in the Theory of Differential and Pseudo-Differential Equations of Parabolic Type.} Birkh\"auser, Basel.

\bibitem{FOT}Fukushima,  M., Oshima, Y., and Takeda, M.  (1994). {\it Dirichlet Forms and Symmetric Markov
Processes.} de Gruyter, Berlin.

\bibitem{GorenfloSurvey} Gorenflo, R. and  Mainardi, F. (2003). Fractional diffusion processes: Probability distribution and continuous time random walk. {\it  Lecture Notes in Physics} {\bf 621} 148--166.

\bibitem{krageloh}  Kr\"{a}geloh, A.M. (2003). Two families of functions related to the fractional powers of generators of strongly continuous contraction semigroups, \emph{J. Math. Anal. Appl.} {\bf 283} 459-467.

\bibitem{fadePRE} Meerschaert, M.M., Benson, D.A., and Baeumer, B. (1999). Multidimensional advection and fractional dispersion. {\it Phys. Rev. E} {\bf 59}, 5026--5028.

\bibitem{RVbook}
 Meerschaert, M.M. and Scheffler, H.-P.  (2001).
\newblock {\it Limit Distributions for Sums of Independent Random Vectors: Heavy Tails in Theory and Practice}.
\newblock Wiley, New York.

\bibitem{Zsolution} Meerschaert, M.M.,  Benson, D.A., Scheffler, H.-P. and Baeumer, B. (2002).
Stochastic solution of space-time fractional diffusion equations.
{\it Phys. Rev. E} {\bf 65}, 1103--1106.

\bibitem{limitCTRW}  Meerschaert, M.M. and  Scheffler, H.-P. (2004).
Limit theorems for continuous time random walks with infinite mean
waiting times. {\it J. Applied Probab.} {\bf 41}, No. 3, 623--638.


\bibitem{M-S-ultra} Meerschaert, M.M. and  Scheffler, H.-P. (2006).  Stochastic model for ultraslow
diffusion. {\it  Stochastic Processes  Appl.} {\bf 116}, 1215--1235.

\bibitem{DOFCP} Meerschaert, M.M., Nane, E., and Vellaisamy, P. (2011). Distributed-order fractional diffusions on bounded domains. {\it J. Math. Anal. Appl.} {\bf 379}, 216--228.


\bibitem{MetzlerKlafter}  Metzler, R. and  Klafter, J. (2004). The restaurant at the end of the random walk: recent developments in the description of anomalous transport by fractional dynamics.  {\it J. Physics A} {\bf 37},  R161--R208.

\bibitem{Podlubny}  Podlubny, I. (1999). { \it Fractional Differential Equations}, Academic Press, San Diego.

\end{thebibliography}
\end{document}